\documentclass[12pt]{amsart}
\usepackage{amsmath, amsthm, amssymb,cite}
\usepackage{fullpage}
\usepackage{color}
\usepackage[colorlinks=true]{hyperref}

\newtheorem{theorem}{Theorem}[section]
\newtheorem{lemma}{Lemma}

\numberwithin{lemma}{section}

\allowdisplaybreaks

\numberwithin{equation}{section}

\author{Jianxiang Liu}
\address{Department of Mathematics \& Statistics, McMaster University}
\curraddr{}
\email{liu2238@mcmaster.ca}

\begin{document}

\title{Interior curvature estimate for curvature quotient equations on convex hypersurfaces}

\begin{abstract}
In this paper, we study interior curvature estimates for convex hypersurfaces satisfying the curvature quotient equation $\frac{\sigma_n}{\sigma_{n-2}}(\kappa) = f(X) > 0$ on Riemannian manifolds. Using a pointwise approach combined with a case-by-case analysis adapted to the structure of the equation, we establish interior $C^2$ estimates in all dimensions $n \geq 3$.
\end{abstract}

\maketitle

\tableofcontents

\section{Introduction}
\label{sec:1}
Interior curvature estimates are known to play a fundamental role in the study of regularity for solutions to geometric partial differential equations. In this paper, we aim to derive interior curvature bounds for convex solutions to a specific class of curvature equations expressed as
\begin{equation}
	\frac{\sigma_n}{\sigma_{n-2}}(\kappa(X)) = f(X), \quad f(X) > 0,
	\label{eq:originaleq}
\end{equation}
where $\sigma_k$ represents the $k$-th elementary symmetric polynomial, and $X,\kappa(X)$ represent the position vector of $M$ and the principal curvatures of hypersurfaces $M\subset \mathbb{R}^{n+1}$ at $X$ respectively. 

For equation (\ref{eq:originaleq}), we establish the following interior $C^2$ estimate. 
\begin{theorem}
For $n\geq3$, suppose that $M=(x,u(x))$ is a local $C^4$ graph over a ball $x\in B_r\subset \mathbb{R}^n$ with positive principal curvatures $\kappa=(\kappa_{1},\cdots,\kappa_{n})>0$ and that it is a solution to equation (\ref{eq:originaleq}) in $B_r$. And $f\in C^{2}(B_r)$ is a positive function. Then there is a positive constant $C$ such that
	$$ \text{sup}_{B_\frac{r}{2}}|\kappa_{i}|\leq C, $$
	where $C$ depends only on $ n,r,\Vert f \Vert_{C^{2}(B_r)}, inf_{B_r}|f|, \Vert M \Vert_{C^{1}(B_r)} $.
	\label{thm:1.1}	
\end{theorem}

Equation (\ref{eq:originaleq}) is a special case of the equations studied in the classical seminal works~\cite{caffarelli1985nonlinear, caffarelli1988nonlinear} by L.~Caffarelli, L.~Nirenberg, and J.~Spruck. They studied a broad class of fully nonlinear operators $\Phi$, which can be represented as smooth symmetric functions of the form $\Phi(D^2 u) = F(\lambda(D^2 u))$, where $\lambda$ denotes the eigenvalues of $D^2 u$.

Our work is more closely related to \cite{caffarelli1988nonlinear}, which focuses on curvature equations. There, it is assumed that $F\in C^2(\Gamma)\cap C^0(\bar{\Gamma})$ is a symmetric function defined on a bounded, convex and symmetric region $\Gamma\subset \mathbb{R}^{n},\Gamma\neq \mathbb{R}^{n}$, with $0\in \partial\Gamma$ and such that $\Gamma+\Gamma_+\subset\Gamma$, where $\Gamma_+$ is the positive cone in $\mathbb{R}^n$ and $u$ is a function defined in $\Gamma$.

At any point $x\in \Gamma$, the principal curvatures $\kappa=(\kappa_1,\cdots,\kappa_n)$ of the graph $X=(x,u(x))$ satisfy curvature equations
\begin{align}
	F[\kappa(X)]=f(X,\nu(X)),
	\label{eq:1.2}
\end{align}
where $ \nu(X) $ is the outer normal vector at $X$.

We are also inspired by ~\cite{caffarelli1985nonlinear}, where a class of fully nonlinear Hessian equations
\begin{align}
	F[\lambda(D^2u)]=g(x,u,Du)
	\label{eq:1.3}
\end{align}
are investigated.

Throughout the literatures, they assume that $F$ satisfies the following structural conditions. The function $F$ is defined in an open convex cone $\Gamma \subset \mathbb{R}^n$ with vertex at the origin and containing the positive cone $\{\lambda \in \mathbb{R}^n \mid \lambda_i > 0 \text{ for all } i\}$, such that
\begin{align}
    F>0 \hspace{3mm} in \hspace{3mm}\Gamma,& \hspace{0.5cm} F=0 \hspace{3mm} on \hspace{3mm} \partial\Gamma,\label{eq:1.4}\\
    F&\text{ is concave in }\Gamma,\label{eq:1.5}\\
    \sum_{i=1}^{n} F^{ii}&\geq \omega_1 \hspace{3mm} on \hspace{3mm} \Gamma_{\mu_1,\mu_2}\label{eq:1.6},\\
    \sum_{i=1}^{n} F^{ii}\lambda_{i}&\geq \omega_2 \hspace{3mm} on \hspace{3mm} \Gamma_{\mu_1,\mu_2}\label{eq:1.7},
\end{align}
where $\Gamma_{\mu_1,\mu_2}=\{\lambda\in\Gamma:\mu_1\leq F[u]\leq \mu_2\}$ for any $\mu_2\geq\mu_1 > 0$, and $\omega_1,\omega_2$ are positive constants depending on $\mu_1,\mu_2$, and $F^{ii}=\frac{\partial F}{\partial\lambda_{i}}$.

As notable examples of this class, the $k$-Hessian operators $\sigma_k^{1/k}$ and the Hessian quotient operators $\left( \frac{\sigma_k}{\sigma_l} \right)^{1/(k-l)}$, where $n \geq k > l \geq 1$, have been extensively studied in geometric analysis due to their close connections with curvature. For instance, $\sigma_k$ corresponds to the mean curvature, scalar curvature, and Gauss curvature when $k = 1$, $2$, and $n$, respectively.

The study of interior curvature and $C^2$ estimates for solutions to Equations~(\ref{eq:1.2}) and~(\ref{eq:1.3}) has been an area of active and ongoing research. This direction is particularly important for understanding regularity of solution in the absence of boundary data. The problem traces back to the work of Heinz~\cite{heinz1959elliptic} in 1959, where he established an interior $C^2$ estimate for the Monge–Ampère equation $\det(D^2 u) = f(x, u, Du)$ in $\mathbb{R}^2$. Subsequently, similar results were obtained by Chen-Han-Ou~\cite{chen2016interior}, and independently by Liu~\cite{Liu2021}, using different methods.

However, Pogorelov demonstrated that the interior $C^2$ estimate does not hold for the Monge--Ampère equation in dimensions $n \geq 3$, via a counterexample~\cite{pogorelov}. Urbas later extended this counterexample to show that the interior $C^2$ estimate also breaks down for the Hessian equation $\sigma_k(D^2 u) = f$ with $k \geq 3$~\cite{urbas1990existence}. Despite these findings, the case of the Hessian equation $\sigma_2(D^2 u) = f$ is not yet fully understood, and its general-dimensional form remains an open problem.

Notable progress has been made on this problem under various structural conditions in dimensions $n \geq 3$. When $n = 3$, Warren and Yuan~\cite{warren2007hessianestimatessigma2equation} established an interior $C^2$ estimate for $\sigma_2(D^2 u) = 1$. Building on this, Qiu extended the result to the more general equation $\sigma_2(D^2 u) = f(x, u, Du)$~\cite{Qiu2018}, and further to the curvature setting $\sigma_2(\kappa) = f(X, \nu(X))$~\cite{qiu2019interiorcurvatureestimateshypersurfaces}. For $n = 4$, an interior estimate for $\sigma_2(D^2 u) = 1$ was obtained by Shankar and Yuan in their recent work~\cite{10.4007/annals.2025.201.2.4}. In higher dimensions, Mooney~\cite{mooney2021strict} established interior estimates for strictly convex solutions of $\sigma_2(D^2 u) = 1$. McGonagle, Song, and Yuan~\cite{mcgonagle2019hessian} further generalized this result to solutions satisfying $D^2 u > -c(n)$. Additional refinements were introduced by Shankar and Yuan under the relaxed assumptions $D^2 u > -K I$~\cite{shankar2020hessian} and $D^2 u > -c(n) \Delta u$~\cite{10.4007/annals.2025.201.2.4}, respectively. A meaningful result was obtained by Guan and Qiu~\cite{10.1215/00127094-2019-0001}, who established interior $C^2$ estimates for both $\sigma_2(D^2 u) = f(x, u, Du)$ and $\sigma_2(\kappa) = f(X, \nu(X))$, under the structural condition $\sigma_3(\kappa) > -A $.

Meanwhile, research on Hessian quotient equations has been relatively limited, due to their complicated structure. For $1 \leq k \leq n$, Sheng–Urbas–Wang~\cite{10.1215/S0012-7094-04-12321-8} and Guan-Zhang~\cite{guanzhangpv1} studied the curvature equations $\frac{\sigma_k}{\sigma_{k-1}}(\kappa) = f(X)$, establishing interior $C^2$ estimates for admissible solutions. For structurally more complicated equations, interior estimates have primarily been established in low dimensions such as three and four. For example, Chen–Warren–Yuan~\cite{chen2009priori} as well as Wang-Yuan~\cite{wang2014hessian}, obtained interior $C^2$ estimates for convex solutions by employing the special Lagrangian structure of the equation $\frac{\sigma_3}{\sigma_1}(D^2 u) = 1$. Zhou~\cite{zhou2024notes} refined these results by incorporating a twisted Lagrangian structure into the equation $\frac{\sigma_3}{\sigma_1}(D^2 u) = f(x)$. The three-dimensional case was addressed by Lu~\cite{lu2023interiorc2estimatehessian} using the Legendre transform. More recently, substantial progress was made in~\cite{lu2024interiorc2estimatehessian}, where Lu extended the interior $C^2$ estimates for convex solutions of equations of the same form to higher dimensions. He also showed that such estimates fail for equations of the type $\frac{\sigma_k}{\sigma_l}(D^2 u) = f$ whenever $1 \leq l < k \leq n$ and $k - l \geq 3$, following the approach of Urbas~\cite{urbas1990existence}. In the latest paper\cite{lu2026liouvillestheoremhessianquotient}, Lu and Sroka solve the interior $C^2$ estimate for $\frac{\sigma_2}{\sigma_1}(D^2 u) =1$, thus the remaining cases for Hessian quotient equations are $\frac{\sigma_{k}}{\sigma_{k-1}}(D^2 u) = f$ and $\frac{\sigma_k}{\sigma_{k-2}}(D^2 u) = f$ with $3 \leq k < n$. Meanwhile, the same estimate problem for general curvature quotient equations $\frac{\sigma_k}{\sigma_l}(\kappa) = f$ for $1 \leq l < k \leq n$ remain open as well.

Inspired by Lu’s work \cite{lu2024interiorc2estimatehessian}, we revisit the Jacobi inequality on Riemannian manifolds. In contrast to the integral method developed by Lu, which draws on the framework of Shankar and Yuan \cite{10.4007/annals.2025.201.2.4}, our approach uses the maximum principle, thus avoiding the technical difficulties of applying the Legendre transform on Riemannian manifolds. Building on the auxiliary function technique of Guan and Qiu \cite{10.1215/00127094-2019-0001}, we develop a new case-by-case analysis tailored to the structure of Hessian quotient equations and obtain a similar result. This fills the gap of the interior curvature estimate for curvature equation $\frac{\sigma_k}{\sigma_{k-2}}(\kappa) = f(X)$ in general dimension.

The organization of the paper is as follows. In Section~\ref{sec:2}, we collect and prove several preliminary results related to the Hessian quotient operator $\frac{\sigma_n}{\sigma_{n-2}}$, including a key concavity inequality. In Section~\ref{sec:3}, we derive the Jacobi inequality for this operator on Riemannian manifolds. Finally, in Section~\ref{sec:4}, we complete the proof of Theorem~\ref{thm:1.1} using the results established in Sections~\ref{sec:2} and~\ref{sec:3}.

\section{Preliminaries} \label{sec:2}
Define the operator $F$ as 
\begin{align*}
	F = \frac{\sigma_{n}}{\sigma_{n-2}}.
\end{align*}
Let $h = (h_{ij})$ denote the second fundamental form of $M$, whose eigenvalues are denoted by $\kappa = (\kappa_1, \dots, \kappa_n) \in \mathbb{R}^n$. We introduce the following notation
\begin{align*}
	(\kappa|i) = (\kappa_1, \dots, \kappa_{i-1}, \kappa_{i+1}, \dots, \kappa_n) \in \mathbb{R}^{n-1},
\end{align*}
where $(\kappa|i)$ corresponds to the vector formed by removing the $i$-th component of $\kappa$. Analogously, $(\kappa|ij)$ denotes the vector obtained by deleting both the $i$-th and $j$-th components of $\kappa$.

We first collect some basic formulas for the Hessian operator (see, for instance, \cite{lin1994dirichlet}).

\begin{lemma}\label{Sigma_k basic}
	For any $\kappa=(\kappa_1,\cdots,\kappa_n)\in \mathbb{R}^n$, the following identities hold:
	\begin{align*}
		\sigma_k(\kappa)= \kappa_i\sigma_{k-1}(\kappa|i)+\sigma_k(\kappa|i),\quad \sum_i\sigma_{k}(\kappa|i)&=(n-k)\sigma_{k}(\kappa),\\
		\sum_i\kappa_i\sigma_{k-1}(\kappa|i)= k\sigma_k(\kappa),\quad \sum_i\kappa_i^2\sigma_{k-1}(\kappa|i)= \sigma_1(\kappa)&\sigma_k(\kappa)-(k+1)\sigma_{k+1}(\kappa).
	\end{align*}
\end{lemma}

We now summarize fundamental properties of the operator $ F $. Proofs of Lemma \ref{F^{ijkl}} was provided in \cite{lu2024interiorc2estimatehessian}.

\begin{lemma}\label{F^{ijkl}}
	Let $1\leq k<n$ and let $n\geq 2$. Assume the coordinates are rotated such that $h$ is diagonalized at $X(x_0)$ for $x_0\in B_r$. Then at $X(x_0)$,
	\begin{align*}
		F^{pq}=\frac{\partial F}{\partial h_{pq}}=\frac{\sigma_n^{pq}}{\sigma_k}-\frac{\sigma_n\sigma_k^{pq}}{\sigma_k^2},
	\end{align*}
	\begin{align*}
		F^{pq,rs}=\frac{\partial^2F}{\partial h_{pq}\partial h_{rs}}= \begin{cases}
			-2\frac{\sigma_n^{pp}\sigma_k^{pp}}{\sigma_k^2}+2\frac{\sigma_n\left(\sigma_k^{pp}\right)^2}{\sigma_k^3},\quad &p=q=r=s,\\
			\frac{\sigma_n^{pp,rr}}{\sigma_k}-\frac{\sigma_n^{pp}\sigma_k^{rr}}{\sigma_k^2}-\frac{\sigma_n^{rr}\sigma_k^{pp}}{\sigma_k^2}-\frac{\sigma_n\sigma_k^{pp,rr}}{\sigma_k^2}+2\frac{\sigma_n\sigma_k^{pp}\sigma_k^{rr}}{\sigma_k^3}, &p=q,r=s,p\neq r,\\
			\frac{\sigma_n^{pq,qp}}{\sigma_k}-\frac{\sigma_n\sigma_k^{pq,qp}}{\sigma_k^2}, & p=s,q=r,p\neq q,\\
			0, & \textit{otherwise}.
		\end{cases}
	\end{align*}
	
	Moreover, for all $p\neq q$, the inequality $-F^{pq,qp}>0$ holds, and
	\begin{align*}
		-F^{pq,qp}=\frac{F^{pp}-F^{qq}}{\kappa_q-\kappa_p},\quad\kappa_p\neq \kappa_q,
	\end{align*}
	where $\kappa_p$ and $\kappa_q$ denote the eigenvalues of 
    $h$. 
\end{lemma}

\begin{lemma}\label{lm:lambda_n-1,nbound}
	Let $n\geq 3$. Suppose $M$ is a convex graph over $B_r\subset\mathbb{R}^n$ with coordinates rotated to diagonalize $h$ at $X(x_0)$, and the eigenvalues arranged as $\kappa_1\geq \cdots\geq \kappa_n>0$. Then at $X(x_0)$,
	\begin{align}\label{eq:lambda_n-1,nbound}
		\kappa_n\leq \sqrt{C(n)F},\quad \kappa_{n-1}\geq \sqrt{\frac{F}{C(n)}},
	\end{align}
	where $C(n)$ is a positive constant depending only on $n$. Furthermore, we also have
    \begin{align}\label{eq:lambda_n-1,nbound+}
		\dfrac{1}{C(n)F}\leq \dfrac{1}{\kappa_{n-1}\kappa_n} \leq \dfrac{1}{F}
	\end{align}    
\end{lemma}
\begin{proof}
The proof of inequality (\ref{eq:lambda_n-1,nbound}) can also be found in \cite{lu2024interiorc2estimatehessian}.

    We just show inequality (\ref{eq:lambda_n-1,nbound+}) here. First, recall the identity for $1\leq l<k\leq n$,
	\begin{align}
		\frac{\sigma_{n-k}}{\sigma_{n-l}}(\kappa)=\dfrac{\sigma_k}{\sigma_{l}}(\dfrac{1}{\kappa})\label{eq:inversesigma}
	\end{align}
where $\kappa^{-1}=(\kappa^{-1}_1,\kappa^{-1}_2,\cdots,\kappa^{-1}_n)$.

Hence, we have
    \begin{align*}
		\frac{1}{F}=\sigma_2(\kappa^{-1})=\dfrac{1}{2}\left(\sum_{i=1}^{n}\sum_{k\neq i}\dfrac{1}{\kappa_{i}\kappa_{k}}\right).
    \end{align*}

In addition, we have the dominant term of $\frac{1}{F}=\sigma_2(\kappa^{-1})$ is $\frac{1}{\kappa_{n-1}\kappa_n}$, since
\begin{align*}
		\dfrac{1}{\kappa_{n-1}\kappa_n}\leq\dfrac{1}{2}\left(\sum_{i=1}^{n}\sum_{k\neq i}\dfrac{1}{\kappa_{i}\kappa_{k}}\right) \leq \dfrac{C(n)}{\kappa_{n-1}\kappa_n}.
    \end{align*}

This completes the proof of inequality (\ref{eq:lambda_n-1,nbound+}).
\end{proof}

\begin{lemma}\label{lm:fibound}
	Under the assumptions of Lemma \ref{lm:lambda_n-1,nbound}, the following inequalities hold at $X(x_0)$:
	\begin{align}
     F^{ii}\geq \dfrac{F^2}{\kappa^2_i\kappa_n} (\text{for }i\neq n);&\quad F^{nn}\geq \dfrac{F^2}{\kappa^2_n\kappa_{n-1}}  \label{eq:fi_ele}\\
		\dfrac{F}{C_1(n)\kappa_{n}}\leq&\sum_{i=1}^{n}F^{ii}\leq \dfrac{2F}{\kappa_{n}}\label{eq:fi}\\
		&\sum_{i=1}^{n}F^{ii}h_{ii}=2F\label{eq:fihi}\\ 
        \dfrac{C_2(n)F^2}{\kappa_{n}}\leq&\sum_{i=1}^{n}F^{ii}h^2_{ii}\leq \dfrac{C_3(n)F^2}{\kappa_{n}}\label{eq:fihi^2}
	\end{align}
	where $C_1(n),C_2(n),C_3(n)$ are positive constants depending only on $n$.
\end{lemma}

\begin{proof}

(\ref{eq:fihi}) follows directly from formula (\ref{eq:inversesigma}). For inequality (\ref{eq:fi_ele}) and (\ref{eq:fi}),
\begin{align*}
    F^{ii}&= \frac{F^2}{\kappa^2_{i}}\cdot\sum_{k\neq i}\dfrac{1}{\kappa_{k}} \\
	\sum_{i=1}^{n}F^{ii}&=\sum_{i=1}^{n}\left( \dfrac{F^2}{\kappa^2_{i}}\cdot\sum_{k\neq i}\dfrac{1}{\kappa_{k}}  \right) \leq \dfrac{F^2}{\kappa_{n}}\sum_{i=1}^{n}\sum_{k\neq i}\dfrac{1}{\kappa_{i}\kappa_{k}} =\dfrac{2F}{\kappa_{n}},\\
	\sum_{i=1}^{n}F^{ii}&=\sum_{i=1}^{n}\left( \dfrac{F^2}{\kappa^2_{i}}\cdot\sum_{k\neq i}\dfrac{1}{\kappa_{k}}  \right) \geq \dfrac{F^2}{\kappa_{n}}\sum_{k\neq n}\dfrac{1}{\kappa_{n}\kappa_{k}}\geq \dfrac{F^2}{n\kappa_{n}}\sum_{i=1}^{n}\sum_{k\neq i}\dfrac{1}{\kappa_{i}\kappa_{k}}=\dfrac{2F}{n\kappa_{n}},
\end{align*}
For inequality (\ref{eq:fihi^2}),
\begin{align*}
	\sum_{i=1}^{n}F^{ii}h^2_{ii}&=\sum_{i=1}^{n}\left( F^2\cdot\sum_{k\neq i}\dfrac{1}{\kappa_{k}}  \right) \leq\sum_{i=1}^{n}\left( F^2\cdot\sum_{k\neq i}\dfrac{1}{\kappa_{n}}  \right)\leq \dfrac{C_3(n)F^2}{\kappa_{n}},\\
	\sum_{i=1}^{n}F^{ii}h^2_{ii}&=\sum_{i=1}^{n}\left( F^2\cdot\sum_{k\neq i}\dfrac{1}{\kappa_{k}}  \right) \geq\sum_{i=1}^{n-1}\left( F^2\cdot\dfrac{1}{\kappa_{n}}  \right)\geq \dfrac{C_2(n)F^2}{\kappa_{n}}.
\end{align*}
This completes the proof.
\end{proof}

From (\ref{eq:inversesigma}), we see that $ 1/(\kappa_{n-1}\kappa_{n})$ is the dominant term in $ \sigma_2(\kappa^{-1})=1/F $. Hence, the following bound to $\kappa_{n}$ will be useful later.

\begin{lemma}\label{lm:lambda_nbound}
	Under the assumptions of Lemma \ref{lm:lambda_n-1,nbound}, at $X(x_0)$,
	\begin{align*}
		\dfrac{1}{\kappa_{1}}+\dfrac{1}{\kappa_{2}}+\cdots+\dfrac{1}{\kappa_{n-1}} \leq \frac{\kappa_n}{F} \leq C(n)\left( \dfrac{1}{\kappa_{1}}+\dfrac{1}{\kappa_{2}}+\cdots+\dfrac{1}{\kappa_{n-1}}\right) 
	\end{align*}
	where $C(n)$ is a positive constant depending only on $n$.
\end{lemma}

\begin{proof}
	By formula (\ref{eq:inversesigma}), we have:
	\begin{align*}
		\frac{1}{F}=\sigma_2(\kappa^{-1})=\dfrac{1}{2}\left(\sum_{i=1}^{n}\sum_{k\neq i}\dfrac{1}{\kappa_{i}\kappa_{k}}\right).
	\end{align*}

    The first inequality is immediate. For the second inequality:
    \begin{align*}
    	\dfrac{2}{F}=2\sigma_2(\kappa^{-1})=\sum_{i=1}^{n-1}\sum_{k\neq i}\dfrac{1}{\kappa_{i}\kappa_{k}}+\sum_{k\neq n}\dfrac{1}{\kappa_{n}\kappa_{k}}\leq\dfrac{C_4(n)}{\kappa_{n-1}\kappa_{n}}+\dfrac{1}{\kappa_{n}}\sum_{k\neq n}\dfrac{1}{\kappa_{k}}\leq\dfrac{C_5(n)}{\kappa_{n}}\sum_{k\neq n}\dfrac{1}{\kappa_{k}}.
    \end{align*}
    where $C_4(n),C_5(n)$ are positive constants depending only on $n$.
\end{proof}

To finish this section, we demonstrate an essential concavity inequality for $F$. This inequality was introduced in \cite{guansroka} first and a proof can be found in \cite{lu2024interiorc2estimatehessian}.

\begin{lemma}\label{lm:concavityineq}
	Let $n\geq 3$ and let $F=\frac{\sigma_n}{\sigma_{n-2}}$ with $\kappa_1\geq\cdots\geq\kappa_n>0$. Then for any $\xi=(\xi_1,\cdots,\xi_n)\in \mathbb{R}^n$, we have
	\begin{align}
		-\sum_{i,j}F^{ii,jj}\xi_i\xi_j-\frac{F^{11}\xi_1^2}{\kappa_1}\geq  -\frac{2}{F}\left(\sum_iF^{ii}\xi_i\right)^2+\frac{1}{2(n-1)}\frac{F^{11}\xi_1^2}{\kappa_1}.\label{eq:concavityineq}
	\end{align}
\end{lemma}

\section{Jacobi inequality}\label{sec:3}
To establish curvature estimates, the primary challenges arise from the first and second derivatives of curvature. In this section, we generalize the Jacobi inequality in \cite{lu2024interiorc2estimatehessian} to Riemannian manifolds using the concavity inequality in Lemma \ref{lm:concavityineq} in order to eliminate the second derivatives of curvature subsequently. 

For simplicity, assume $\kappa_{1}$ is sufficiently large.
Let $M\subset\mathbb{R}^{n+1}$ be a hypersurface represented as a graph $X=(x,u(x))$ over $B_r\subseteq\mathbb{R}^n$ satisfying the equation (\ref{eq:originaleq}). At any point $x\in B_r$, let $h=(h_{ij})$ denote the second fundamental form of the graph and  $\kappa=(\kappa_{1},\kappa_{2},\cdots,\kappa_{n})$ be its principal curvatures. 

We now recall several fundamental identities in Riemannian geometry. First, we choose an orthonormal frame in $\mathbb{R}^{n+1}$ such that $\{ e_1,e_2,\cdots,e_n \}$ are tangent to $M$ and $\nu$ is the outer normal to $M$. The following fundamental identities hold:
\begin{align}
	X_{ij}&=-h_{ij}\nu \quad \text{(Gauss formula)} ,\label{eq:Gauss formula}\\
	\nu_i&=\sum_{j=1}^{n}h_{ij}e_j \quad \text{(Weingarten equation)}, \label{eq:Weingarten equation}\\
	h_{ijk}&=h_{ikj} \quad \text{(Codazzi equation)}, \label{eq:Codazzi equation}\\
	R_{ijkl}&=h_{ik}h_{jl}-h_{il}h_{jk} \quad \text{(Gauss equation)}\label{eq:Gauss equation},
\end{align}
where $R_{ijkl}$ is the curvature tensor. Define $h_{ijkl}=\nabla_l\nabla_kh_{ij}$, the commutator formula yields:
\begin{align}
	h_{ijkl}-h_{ijlk}=\sum_{m}h_{mi}R_{mjkl}+\sum_{m}h_{mj}R_{mikl}.\label{eq:commutator}  
\end{align}
Combine (\ref{eq:commutator}) with Codazzi and Gauss equation, we derive:
\begin{align}
	h_{iikk}&=h_{kkii}+\sum_{m}(h_{im}h_{mi}h_{kk}-h^2_{mk}h_{ii}).  \label{eq:commutation}
\end{align}

We now proceed to generalize the Jacobi inequality.

\begin{lemma}\label{lm:Jacobi}
	For $n\geq 3$, suppose $f\in C^{1,1}(B_{r})$ be a positive function and let $M$ be a convex graph over $B_r\subset \mathbb{R}^n$ solving equation (\ref{eq:originaleq}). At any $X(x_0),x_0\in B_{r}$, assume $h$ is diagonalized with ordered principal curvatures $\kappa_1\geq \cdots\geq \kappa_n$. Define $b=\ln\kappa_1$. Then, at $X(x_0)$, the following inequality holds in the viscosity sense:
	\begin{align}
		\sum_i F^{ii}b_{ii}\geq \hat{c}(n)\sum_iF^{ii}b_i^2+\sum_{i}F^{ii}h_{ii}h_{11}-\sum_{i}F^{ii}h^{2}_{ii}-C,\label{eq:Jacobi}
	\end{align}
	where $\hat{c}(n)$ is a positive constant depending only on $n$, and $C$ depends only on $n$, $\|M\|_{C^{0,1}\left( B_r\right) }$, $\min_{B_r}f $, and $\|f\|_{C^{1,1}\left( B_r\right)}$.
\end{lemma} 

\begin{proof}
	
	Suppose $\kappa_1$ has multiplicity $m$ at $X(x_0)$. By Lemma 5 in \cite{BCD}, we have
	\begin{align}
		\delta_{\alpha\beta}\cdot{(\kappa_1)}_i=h_{\alpha\beta i},\quad 1\leq \alpha,\beta\leq m,\label{BCD1}\\
		{(\kappa_1)}_{ii}\geq h_{11ii}+2\sum_{p>m}\frac{h_{1pi}^2}{\kappa_1 -\kappa_p},\label{BCD2}
	\end{align}
	in the viscosity sense.
	
	It follows that
	\begin{align*}
		b_{ii}=\frac{{(\kappa_1)}_{ii}}{\kappa_1}-\frac{{(\kappa_1)}_i^2}{\kappa_1^2}\geq \frac{h_{11ii}}{ \kappa_1}+2\sum_{p>m}\frac{h_{1pi}^2}{\kappa_1 (\kappa_1 -\kappa_p )}-\frac{ h_{11i}^2}{\kappa_{1}^2}.
	\end{align*}
	
	Contracting with $F^{ii}$,
	\begin{align}\label{eq:fibii-0}
		\sum_i F^{ii}b_{ii}\geq \sum_i \frac{F^{ii}h_{11ii}}{ \kappa_1}+2\sum_i\sum_{p>m}\frac{F^{ii}h_{1pi}^2}{\kappa_1 (\kappa_1 -\kappa_p )}-\sum_i \frac{F^{ii} h_{11i}^2}{\kappa_{1}^2}.
	\end{align}
	
	Differentiate equation (\ref{eq:originaleq}) twice with respect to $e_1$:
	\begin{align}
		F^{ij}h_{ij1}&=d_Xf(e_1)\label{eq:1deriofeq},\\
		F^{ij}(h_{ij11})+F^{pq,rs}h_{pq1}h_{rs1}&=d_{X,X}f(e_1,e_1)\geq-C\kappa_{1}-C,\label{eq:2deriofeq}
	\end{align}
    where $C$ is a positive constant depending only on $n$, $\|M\|_{C^{0,1}\left( B_r\right) }$, $\min_{B_r}f $, and $\|f\|_{C^{1,1}\left( B_r\right)}$ and may change from line to line.
    
    Substituting (\ref{eq:commutation}) and (\ref{eq:2deriofeq}) into (\ref{eq:fibii-0}),
    \begin{align}\label{eq:fibii-1}
    	\sum_i F^{ii}b_{ii}\geq &-\sum_{p,q,r,s} \frac{F^{pq,rs}h_{pq1}h_{rs1}}{\kappa_1 }+2\sum_i\sum_{p>m}\frac{F^{ii}h_{1pi}^2}{\kappa_1 (\kappa_1 -\kappa_p )}-\sum_i \frac{ F^{ii}h_{11i}^2}{\kappa_{1}^2}\notag\\
    	&+\sum_{i}F^{ii}h_{ii}h_{11}-\sum_{i}F^{ii}h^{2}_{ii}-C.
    \end{align}

    By Lemma \ref{F^{ijkl}},
    \begin{align*}
    	-\sum_{p,q,r,s} F^{pq,rs}h_{pq1}h_{rs1}=-\sum_{i,j}F^{ii,jj}h_{ii1}u_{jj1}-\sum_{i\neq j}F^{ij,ji}h_{ij1}^2.
    \end{align*}
    
    Substituting into (\ref{eq:fibii-1}),
    \begin{align}\label{eq:fibii-2}
    	\sum_i F^{ii}b_{ii}\geq &\ -\sum_{i, j} \frac{F^{ii,jj}h_{ii1}h_{jj1}}{\kappa_1} -\sum_{i\neq j} \frac{F^{ij,ji}h_{ij1}^2}{\kappa_1}+2\sum_i\sum_{p>m}\frac{F^{ii}h_{1pi}^2}{\kappa_1 (\kappa_1 -\kappa_p )}\\\nonumber
    	&\ -\sum_i \frac{ F^{ii}h_{11i}^2}{\kappa_{1}^2}+\sum_{i}F^{ii}h_{ii}h_{11}-\sum_{i}F^{ii}h^{2}_{ii}-C.
    \end{align}
    
    By Lemma \ref{F^{ijkl}},
    \begin{align*}
    	-\sum_{i\neq j} \frac{F^{ij,ji}h_{ij1}^2}{\kappa_1}\geq -2\sum_{i>m} \frac{F^{i1,1i} h_{i11}^2}{\kappa_1}=2\sum_{i>m} \frac{(F^{ii}-F^{11})h_{11i}^2}{\kappa_1(\kappa_1-\kappa_i)}.
    \end{align*}
    
    On the other hand,
    \begin{align*}
    	2\sum_i\sum_{p>m}\frac{F^{ii}h_{1pi}^2}{\kappa_1 (\kappa_1 -\kappa_p )} \geq  2\sum_{p>m}\frac{F^{11}h_{1p1}^2}{\kappa_1 (\kappa_1 -\kappa_p )}=2\sum_{i>m}\frac{F^{11}h_{11i}^2}{\kappa_1 (\kappa_1 -\kappa_i )}.
    \end{align*}
    
    Combining the above two inequalities and substituting into inequality (\ref{eq:fibii-2}),
    \begin{align}\label{eq:fibii-3}
    	\sum_i F^{ii}b_{ii}\geq &\ -\sum_{i,j} \frac{F^{ii,jj} h_{ii1}h_{jj1}}{\kappa_1}+2\sum_{i>m}\frac{F^{ii}h_{11i}^2}{\kappa_1 (\kappa_1 -\kappa_i )}-\sum_i \frac{ F^{ii}h_{11i}^2}{\kappa_{1}^2}\notag\\
    	&+\sum_{i}F^{ii}h_{ii}h_{11}-\sum_{i}F^{ii}h^{2}_{ii}-C.
    \end{align}
    
    By (\ref{BCD1}) and Codazzi equation,
    \begin{align}\label{BCD1'}
    	h_{11i}=h_{1i1}=\delta_{1i}\cdot(\kappa_1)_1=0,\quad 1<i\leq m.
    \end{align}
    
    Substituting into inequality (\ref{eq:fibii-3}),
    \begin{align*}
    	\sum_i F^{ii}b_{ii}\geq  &-\sum_{i,j} \frac{F^{ii,jj} h_{ii1}h_{jj1}}{\kappa_1}+\sum_{i>m}\frac{F^{ii}h_{11i}^2}{\kappa_1^2}- \frac{ F^{11}h_{111}^2}{\kappa_{1}^2}\notag\\
    	&+\sum_{i}F^{ii}h_{ii}h_{11}-\sum_{i}F^{ii}h^{2}_{ii}-C.
    \end{align*}
    
    Together with Lemma \ref{lm:concavityineq} and (\ref{BCD1'}), then:
    \begin{align*}
    	\sum_i F^{ii}b_{ii}\geq & -\frac{2}{\kappa_1F}\left(\sum_iF^{ii}h_{ii1}\right)^2+\sum_{i>m}\frac{F^{ii}h_{11i}^2}{\kappa_1^2}+\hat{c}(n)\frac{ F^{11}h_{111}^2}{\kappa_{1}^2}\\
    	&+\sum_{i}F^{ii}h_{ii}h_{11}-\sum_{i}F^{ii}h^{2}_{ii}-C\\
    	\geq &\ -\frac{2f_1^2}{\kappa_1F}+\hat{c}(n)\sum_i \frac{F^{ii}h_{11i}^2}{\kappa_1^2}+\sum_{i}F^{ii}h_{ii}h_{11}-\sum_{i}F^{ii}h^{2}_{ii}-C\\
    	\geq &\ \hat{c}(n)\sum_iF^{ii}b_i^2+\sum_{i}F^{ii}h_{ii}h_{11}-\sum_{i}F^{ii}h^{2}_{ii}-C.
    \end{align*}
    The lemma is now proved.

\end{proof}

\section{Proof of the Theorem \ref{thm:1.1}} \label{sec:4}
In this section, we complete the proof of Theorem \ref{thm:1.1} using an appropriate auxiliary function and subsequent calculations. Throughout this section, we adopt the Einstein summation convention unless otherwise stated.

\begin{proof}
    Without loss of generality, we assume $r=1$. The general case $r>0$ follows by a standard rescaling argument. Following the assumption in Section \ref{sec:3}, we consider the same auxiliary function in $B_{1}$ as Guan, Qiu's \cite{10.1215/00127094-2019-0001}, which is:
	\begin{align}\label{testfct}
		P(X(x))=2\log\rho(X)+\log\log \kappa_1-\beta\frac{(X,\nu)}{(\nu,E_{n+1})}+\alpha\frac{1}{(\nu,E_{n+1})^{2}}.
	\end{align}
    Here $E_{n+1}=(0,\cdots,0, 1)$ and $\rho(X)=1-|X|_{\mathbb{R}^{n+1}}^{2}+(X,E_{n+1})^{2}=1-|x|_{\mathbb{R}^{n}}^{2}$.
    The positive constants $\alpha,\beta$ will be determined later. Consequently, the maximum
    of $P$ is attained at an interior point $x_{0}$ of $B_{1}$.
    
    From now on, all the calculation is on $x_{0}$. Since $M$ is a local graph, direct calculation shows
    \begin{align*}
	  |(\nu,E_{n+1})|=\left|\dfrac{\left<(-u_1,\cdots,-u_n,1), E_{n+1} \right>}{\sqrt{1+|Du|^2}}\right|=\dfrac{1}{\sqrt{1+|Du|^2}},
    \end{align*}
    where $\left< \cdot,\cdot \right>$ is the standard inner product in $\mathbb{R}^{n+1}$. Therefore, $|(\nu,E_{n+1})|$ is bounded above and below whenever $\|M\|_{C^{0,1}\left( B_1\right) }$ is bounded.
    
    Using the bounds for $|(\nu,E_{n+1})|$ and assuming that $\kappa_{1}$ is sufficiently large, we obtain:
	\begin{align*}
		C &\leq  \rho^2(x_0)\cdot \log \kappa_1,
	\end{align*}
	that is,
	\begin{align}
		\rho^2 &\geq \frac{C}{\log \kappa_1}, \quad 
		\rho \geq \frac{C}{ \sqrt{\log \kappa_1}}.\label{eq:relaofrhoandlambda_1}
	\end{align}
    where $C$ is a positive constant dependent only on $n$, $|M|_{C^{1}\left( B_{1}\right) }$, $\min_{B_{1}}|f| $, and $\|f\|_{C^{2}\left( B_{1} \right)}$, which may change between steps.
	
    Differentiating $P(X)$ at $x_0 \in B_{1}$, we obtain the first derivative:
	\begin{align}
		0 &= P_i = \frac{2 \rho_i}{\rho} + \frac{h_{11i}}{h_{11} \log h_{11}} 
		- \beta \left[ \frac{(X, \nu)}{(\nu, E_{n+1})}\right] _i 
		- \alpha \frac{2 (\nu_i, E_{n+1})}{(\nu, E_{n+1})^3}. \label{eq:system1}
	\end{align}
    Since $h$ is diagonalized at $x_0$, the second derivative satisfies:
	\begin{align}
		0 \geq F^{ij} P_{ij} = &\frac{2F^{ii} \rho_{ii}}{\rho} 
		- \frac{2 F^{ii} \rho_{i}^2}{\rho^2} 
		- \beta F^{ii} \left[ \frac{(X, \nu)}{(\nu, E_{n+1})}\right] _{ii}+ \frac{F^{ii} (\kappa_{1})_{ii}}{h_{11} \log h_{11}}\notag\\
		& -\left( \log\kappa_{1}+1 \right)  \frac{F^{ii} h^2_{11i}}{h^2_{11} \log^2 h_{11}} 
		- \alpha \frac{2 F^{ii} (\nu_{ii}, E_{n+1})}{(\nu, E_{n+1})^3} 
		+ 6 \alpha \frac{F^{ii} (\nu_i, E_{n+1})^2}{(\nu, E_{n+1})^4}.\label{eq:system2}
	\end{align}

As noted in Section \ref{sec:3}, the critical barrier in the estimate lies in controlling the derivative of curvature. To manage these terms we express certain elements in inequality (\ref{eq:system2}) using derivatives of curvature:
\begin{align}
	\nu_{ii} &= \left( \sum_k h_{ik} e_k \right)_i = \sum_k h_{iik} e_k - \sum_k h_{ik} h_{ki} \nu,\label{eq:systemtm1} \\ \label{eq:systemtm2}
	\left[ \frac{(X, \nu)}{(\nu, E_{n+1})}\right] _i &= \sum_l \frac{(X, e_l) h_{il}}{(\nu, E_{n+1})} 
	- \sum_l \frac{(e_l, E_{n+1}) h_{il} (X, \nu)}{(\nu, E_{n+1})^2}, \\\label{eq:systemtm3}
	\left[ \frac{(X, \nu)}{(\nu, E_{n+1})}\right]_{ii}  &= \frac{h_{ii}}{(\nu, E_{n+1})} 
	- 2\left[  \frac{(X, \nu)}{(\nu, E_{n+1})}\right] _i\cdot \sum_l \frac{(e_l, E_{n+1}) h_{il}}{(\nu, E_{n+1})} \notag\\
	&\quad + \sum_l \frac{(X, e_l) (\nu, E_{n+1}) - (e_l, E_{n+1}) (X, \nu)}{(\nu, E_{n+1})^2} h_{iil}, \\\label{eq:systemtm4}
	\rho_i &= -2 (X, e_i) + 2 (X, E_{n+1}) (e_i, E_{n+1}), \\
	\rho_{ii} &= -2 + 2 (X, \nu) h_{ii} + 2 (e_i, E_{n+1})^2 - 2 (X, E_{n+1}) h_{ii} (\nu, E_{n+1}),\label{eq:systemtm5}
\end{align}

Insert (\ref{eq:systemtm1}---\ref{eq:systemtm5}) and (\ref{eq:system1}) into (\ref{eq:system2}), we obtain:
\begin{align}
    F^{ij} P_{ij} \geq
    &-\dfrac{4\sum_{i}F^{ii}}{\rho}+ \frac{4 F^{ii} h_{ii}[(X, \nu) - (X, E_{n+1}) (\nu, E_{n+1})]}{\rho}+\dfrac{4\sum_iF^{ii}(e_i,E_{n+1})^2}{\rho} \notag\\
    & + \frac{F^{ii} (\kappa_{1})_{ii}}{h_{11} \log h_{11}} 
    -\frac{8F^{ii}[(X, e_i) - (X, E_{n+1}) (e_i, E_{n+1})]^2}{\rho^2} \notag\\
    & - (\log \kappa_{1} + 1) \frac{F^{ii} h^2_{11i}}{h_{11}^2 \log^2 h_{11}} 
    + 2F^{ii}\dfrac{(e_i,E_{n+1})h_{ii}}{(\nu,E_{n+1})}\left( \dfrac{2\rho_{i}}{\rho}+\frac{h_{11i}}{h_{11} \log h_{11}}\right)  \notag\\
    & -\alpha \sum_{k=1}^{n}\frac{2F^{ii} h_{iik} (e_k, E_{n+1})}{(\nu, E_{n+1})^3} - \beta \dfrac{F^{ii} h_{ii}}{(\nu,E_{n+1})} +2\alpha\sum_i
    \frac{F^{ii}h_{ii}^2(e_i,E_{n+1})^2}{(\nu,E_{n+1})^4}\notag\\
    & - \beta F^{ii}h_{iil} \frac{(X,e_l)(\nu, E_{n+1})- (e_l,E_{n+1})(X, \nu)}{(\nu, E_{n+1})^2} 
    +\dfrac{2\alpha F^{ii} h^2_{ii}}{(\nu, E_{n+1})^2},\label{eq:system2-1}
\end{align}

Next, we estimate some terms in inequality (\ref{eq:system2-1}). By  (\ref{eq:fihi}), we have:

\begin{align}
	\frac{4 F^{ii} h_{ii}[(X, \nu) - (X, E_{n+1}) (\nu, E_{n+1})]}{\rho}=& \dfrac{8F[(X, \nu) - (X, E_{n+1}) (\nu, E_{n+1})]}{\rho}\notag\\
	\geq& -\dfrac{C(\parallel f \parallel_{C^2},\parallel M \parallel_{C^{0,1}})}{\rho}.\label{eq:system2-1tm1}
\end{align}

    We apply the Jacobi inequality to eliminate the curvature second-derivative term yields
	\begin{align}
		\dfrac{F^{ii} (\kappa_1)_{ii}}{h_{11}\log h_{11}} 
		&= \frac{F^{ii} (\log \kappa_1)_{ii}}{\log h_{11}} + \frac{1}{\log h_{11}} \cdot \frac{F^{ii} (\kappa_1)_{i}^2}{\kappa_1^2} \notag\\
		&\geq \frac{1}{\log h_{11}} \left( \hat{c}(n) F^{ii} (\log \kappa_1)_i^2 + 2F h_{11} 
		- \sum_{i} F^{ii} h_{ii}^2 - C \right) +\frac{1}{\log h_{11}} \cdot \frac{F^{ii} (\kappa_1)_{i}^2}{\kappa_1^2}\notag\\
		&= \frac{\hat{c}(n)+1}{\log h_{11}}\cdot \frac{F^{ii} (\kappa_1)_{i}^2}{\kappa_1^2}
		+ \frac{2F h_{11}}{\log h_{11}} 
		- \frac{\sum_{i} F^{ii} h_{ii}^2}{\log h_{11}} - C,\label{eq:system2-1tm2}
	\end{align}
    where the constant $\hat{c}(n)$ in this formula comes from Jacobi inequality.

    Combine (\ref{BCD1}), inequality (\ref{eq:system2-1}---\ref{eq:system2-1tm2}) and denote $a_i = [(X, e_i) - (X, E_{n+1})(e_i, E_{n+1})]$, we obtain:
	\begin{align}
		F^{ij} P_{ij} &\geq -\dfrac{4 \sum_{i} F^{ii}}{\rho} - \frac{C}{\rho} - C - \frac{8 \sum_{i}F^{ii} \left[ (X, e_i) - (X, E_{n+1})(e_i, E_{n+1}) \right]^2}{\rho^2} \notag \\
		&\quad + \frac{\hat{c}(n)+1}{\log h_{11}} \cdot \frac{F^{ii} (\kappa_1)_{i}^2}{\kappa_1^2} 
		+ \frac{2F h_{11}}{\log h_{11}} - \frac{\sum_{i} F^{ii} h_{ii}^2}{\log h_{11}} - (\log h_{11} + 1)\dfrac{F^{ii}h^2_{11i}}{h^2_{11} \log^2 h_{11}} \notag \\
		&\quad + \frac{2F^{ii}h_{ii}(e_i, E_{n+1})}{(\nu, E_{n+1})} \left( \frac{2 \rho_{i}}{\rho} 
		+ \frac{h_{11i}}{h_{11} \log h_{11}} \right)
		+ 2 \alpha\frac{ F^{ii} h_{ii}^2}{(\nu, E_{n+1})^2} - C\beta - C\alpha. \notag \\
		&\geq -\dfrac{4 \sum_{i} F^{ii}}{\rho} - C(\frac{1}{\rho} + \beta + \alpha ) -\dfrac{8F^{ii}a_i^2}{\rho^2}
		+ (\hat{c}(n)\cdot\log h_{11}-1)\dfrac{F^{ii}h^2_{11i}}{h^2_{11} \log^2 h_{11}}  \notag \\
		&\quad +\frac{2F h_{11}}{\log h_{11}} - \frac{F^{ii} h_{ii}^2}{\log h_{11}} 
		 + \frac{2F^{ii}h_{ii}(e_i, E_{n+1})}{(\nu, E_{n+1})} \left( \frac{2 \rho_{i}}{\rho} 
		+ \frac{h_{11i}}{h_{11} \log h_{11}} \right) \notag \\
		&\quad + 2 \alpha \frac{F^{ii} h_{ii}^2}{(\nu, E_{n+1})^2}.\label{eq:system2-2}
	\end{align}

By the Cauchy-Schwarz inequality, we have
\begin{align}
	4F^{ii}\frac{h_{ii}\rho_{i}}{\rho(\nu,E_{n+1})}&\geq-\frac{2F^{ii}h_{ii}^{2}}{(\nu,E_{n+1})^{2}}-\frac{2F^{ii}\rho_{i}^{2}}{\rho^{2}},\label{eq:CS1}\\
	2F^{ii}\frac{h_{11i}h_{ii}}{(\nu,E_{n+1})h_{11}\log h_{11}}&\geq-\frac{F^{ii}h_{ii}^{2}}{(\nu,E_{n+1})^{2}}-\frac{F^{ii}h_{11i}^{2}}{h_{11}^{2}\log^{2}h_{11}}.\label{eq:CS2}
\end{align}
	
Thus, substitute inequality (\ref{eq:CS1}) and (\ref{eq:CS2}) into (\ref{eq:system2-2}), we attain
\begin{align}
F^{ij} P_{ij} &\geq -\dfrac{4 \sum_{i} F^{ii}}{\rho} -\dfrac{16F^{ii}a_i^2}{\rho^2}- C(\frac{1}{\rho} + \beta + \alpha ) + (\hat{c}(n)\cdot\log h_{11}-2)\dfrac{F^{ii}h^2_{11i}}{h^2_{11} \log^2 h_{11}}  \notag \\
&\quad +\frac{2F h_{11}}{\log h_{11}} - \frac{F^{ii} h_{ii}^2}{\log h_{11}} 
+ (2 \alpha-3) \frac{F^{ii} h_{ii}^2}{(\nu, E_{n+1})^2}.\label{eq:system2-3}
\end{align}

Assuming that $h_{11}$ is sufficiently large and choosing $\alpha$ sufficiently large, we apply inequality (\ref{BCD2}) to derive
\begin{align}
	F^{ij} P_{ij} &\geq -\dfrac{4 \sum_{i} F^{ii}}{\rho} -\dfrac{16F^{ii}a_i^2}{\rho^2}- C(\frac{1}{\rho} + \beta + \alpha ) 
	+ \dfrac{\hat{c}(n) \log h_{11}}{2}\cdot \dfrac{F^{ii} h_{11i}^2}{h_{11}^2\log^2 h_{11}}  \notag \\
	&\quad+\frac{2F h_{11}}{\log h_{11}} +\alpha \frac{F^{ii} h_{ii}^2}{(\nu, E_{n+1})^2}.\label{eq:system2-4}
\end{align}

We then expand the curvature derivative term in \eqref{eq:system2-4} by the first derivative identity (\ref{eq:system1}). This allows us to absorb one of the negative terms in this inequality, which yields

\begin{align}
	0\geq F^{ij} P_{ij} &\geq -\dfrac{4 \sum_{i} F^{ii}}{\rho} -\dfrac{16F^{ii}a_i^2}{\rho^2}- C(\frac{1}{\rho} + \beta + \alpha ) +\frac{2F h_{11}}{\log h_{11}} +\alpha \frac{F^{ii} h_{ii}^2}{(\nu, E_{n+1})^2} \notag \\
	&\quad+\frac{8\hat{c}(n)\log h_{11}}{\rho^2}\cdot F^{ii}a_i^2 
                \notag \\
&\quad+
\frac{\hat{c}(n)\log h_{11}\cdot F^{ii}h_{ii}^2}{2}
\left(
\beta\frac{b_i}{(\nu,E_{n+1})^2}
+
2\alpha\frac{(e_i,E_{n+1})}{(\nu,E_{n+1})^3}
\right)^2 \notag \\
&\quad-
\frac{4\hat{c}(n)\log h_{11}\cdot F^{ii}a_i h_{ii}}{\rho}
\left|
\beta\frac{b_i}{(\nu,E_{n+1})^2}
+
2\alpha\frac{(e_i,E_{n+1})}{(\nu,E_{n+1})^3}
\right|,
\label{eq:system2-5}
\end{align}

where $b_{i} = (X, e_{i})(\nu, E_{n+1}) - (E_{n+1}, e_{i})(X, \nu)$.

Using (\ref{eq:fihi}) together with (\ref{eq:relaofrhoandlambda_1}), by the assumption that $\kappa_1$ is sufficiently large, we derive

\begin{align}
	0\geq F^{ij} P_{ij} &\geq -\dfrac{4 \sum_{i} F^{ii}}{\rho} -\dfrac{16F^{ii}a_i^2}{\rho^2}- C(\frac{1}{\rho} + \beta + \alpha ) +\frac{2F h_{11}}{\log h_{11}} +\alpha \frac{F^{ii} h_{ii}^2}{(\nu, E_{n+1})^2} \notag \\
	&\quad+\frac{8\hat{c}(n)\log h_{11}}{\rho^2}\cdot F^{ii}a_i^2 
                \notag \\
&\quad+
\frac{\hat{c}(n)\log h_{11}\cdot F^{ii}h_{ii}^2}{2}
\left(
\beta\frac{b_i}{(\nu,E_{n+1})^2}
+
2\alpha\frac{(e_i,E_{n+1})}{(\nu,E_{n+1})^3}
\right)^2 \notag \\
&\quad-
C(n,\alpha,\beta,F,\left\| M \right\|_{C^1})\log^{3/2} h_{11} \notag \\
0&\geq -\dfrac{4 \sum_{i} F^{ii}}{\rho} - C(\frac{1}{\rho} + \beta + \alpha ) +\frac{F h_{11}}{\log h_{11}} +\alpha \frac{F^{ii} h_{ii}^2}{(\nu, E_{n+1})^2} \label{eq:system2-final}\\
	&\quad+\frac{4\hat{c}(n)\log h_{11}}{\rho^2}\cdot F^{ii}a_i^2 \notag \\
&\quad+
\frac{\hat{c}(n)\log h_{11}\cdot F^{ii}h_{ii}^2}{2}
\left(
\beta\frac{b_i}{(\nu,E_{n+1})^2}
+
2\alpha\frac{(e_i,E_{n+1})}{(\nu,E_{n+1})^3}
\right)^2.\notag
\end{align}

Thus, we arrive at the required inequality, in which no negative terms of curvature-derivative remain. From (\ref{eq:system2-final}), we observe that the only dominant negative term is the first one. The consequent discussion on $ a_{i}, b_i $ and $\rho$ provides a range of these two negative terms which will be useful later.
	
To simplify the remaining computations, we now switch to the ambient orthonormal coordinate system $\{E_{1},E_{2},\cdots,E_{n},E_{n+1}\}$
in $\mathbb{R}^{n+1}$, where $E_{i}\perp E_{n+1}$ for all $ i\in \{1,\cdots,n\}$. In this coordinate system, the position vector $X$ admits the decomposition
\begin{equation}
	X=\sum_{i=1}^{n}(X,E_{i})E_{i}+(X,E_{n+1})E_{n+1}.
\end{equation}
Thus,
\begin{equation}\label{eq:neqdefrho}
	\rho=1-\sum_{i=1}^{n}(X,E_{i})^{2}.
\end{equation}
Similarly, we have the decomposition
\begin{align}
	(\nu,E_{n+1})^2+\sum_{i=1}^{n}(e_i,E_{n+1})^2 &=1,\label{eq:sumofE_n+1}\\
    (\nu,E_{n+1})\,\nu+\sum_{i=1}^{n}(e_i,E_{n+1})\,e_i &=E_{n+1} \label{eq:decomposeE_n+1}.
\end{align}

Recalling that
\begin{align*}
	a_{j} &=(X,e_{j})-(X,E_{n+1})(e_{j},E_{n+1})=\sum_{i=1}^{n}(X,E_{i})(E_{i},e_{j}),\\
    b_{i} &= (X, e_{i})(\nu, E_{n+1}) - (E_{n+1}, e_{i})(X, \nu).
\end{align*}

Using (\ref{eq:sumofE_n+1}) and (\ref{eq:decomposeE_n+1}), a direct computation yields
\begin{align*}
\sum_{i=1}^{n}a_i(e_i,E_{n+1})
&=\sum_{i=1}^{n}\left[(X,e_i)(e_i,E_{n+1})
-(X,E_{n+1})(e_i,E_{n+1})^2\right] \\
&=\left(X,\sum_{i=1}^{n}(e_i,E_{n+1})e_i\right)
-(X,E_{n+1})\left(1-(\nu,E_{n+1})^2\right) \\
&=-(\nu,E_{n+1})
\left[(X,\nu)-(X,E_{n+1})(\nu,E_{n+1})\right].
\end{align*}

Consequently,
\begin{align*}
b_i
&=a_i(\nu,E_{n+1})
-(e_i,E_{n+1})
\left[(X,\nu)-(X,E_{n+1})(\nu,E_{n+1})\right] \\
&=a_i(\nu,E_{n+1})
+\frac{(e_i,E_{n+1})}{(\nu,E_{n+1})}
\sum_{j=1}^{n}a_j(e_j,E_{n+1}).
\end{align*}

Equivalently, this relation admits the matrix form
\begin{equation}\label{eq:b=Ma}
    \textbf{b}=T\textbf{a},
\end{equation}
where $\textbf{a} =(a_1,a_2,\cdots,a_n)^t, \textbf{b} =(b_1,b_2,\cdots,b_n)^t$ and if denoting $\textbf{c} =((e_1,E_{n+1}),\cdots,(e_n,E_{n+1}))^t$ we can represent $T$ as:
\begin{align}\label{eq:M}
T=(\nu,E_{n+1})I+\dfrac{\textbf{c}\textbf{c}^t}{(\nu,E_{n+1})}.
\end{align}

One readily verifies that $T$ is a positive-definite, symmetric and also invertible matrix with eigenvalues as 
\begin{align}
\left\{
\underbrace{(\nu,E_{n+1}),\ldots,(\nu,E_{n+1})}_{n-1},
(\nu,E_{n+1})^{-1}
\right\}.
\end{align}

Consequently, $\textbf{b}=T\textbf{a}$ immediately implies

\begin{align}\label{eq:modulea=b}
(\nu,E_{n+1})^2\sum_{i=1}^{n}a_{i}^{2}
\le
\sum_{i=1}^{n}b_{i}^{2}
\le
\frac{\sum_{i=1}^{n}a_{i}^{2}}{(\nu,E_{n+1})^2}.
\end{align}

Using the definition of $a_j$,
\begin{eqnarray*}
	\sum_{i=1}^{n}a_{i}^{2} & = & \sum_{i,k,l=1}^{n}(X,E_{k})(X,E_{l})(E_{k},e_{i})(E_{l},e_{i})\\
	& = & \sum_{k,l=1}^{n}(X,E_{k})(X,E_{l})(\delta_{kl}-(E_{k},\nu)(E_{l},\nu)).
\end{eqnarray*}

Therefore,
\begin{equation}\label{eq:moequvalent-a}
	\frac{1}{C_6}\sum_{i=1}^{n}(X,E_{i})^{2}\leq\sum_{i=1}^{n}a_{i}^{2}\leq C_6\sum_{i=1}^{n}(X,E_{i})^{2},
\end{equation}
where $C_6$ is a positive constant only depend on the upper and lower bound of $(\nu,E_{n+1})$.

Furthermore, inequality (\ref{eq:modulea=b}) shows that $|\textbf{b}|$ also has the similar bound:
\begin{equation}\label{eq:moequvalent-b}
	\frac{1}{C_7}\sum_{i=1}^{n}(X,E_{i})^{2}\leq\sum_{i=1}^{n}b_{i}^{2}\leq C_7\sum_{i=1}^{n}(X,E_{i})^{2}.
\end{equation}
where $C_7$ is a positive constant only depend on the upper and lower bound of $(\nu,E_{n+1})$ as well.

Therefore, our proof strategy is to use positive terms in inequality (\ref{eq:system2-final}) to control the dominant negative term. We categorize the proof into three cases based on the values of $ \rho $ and $b_{i}$. 

~\\
\noindent \textit{\textbf{Case 1}}: \textbf{$\sum_{i=1}^n |(X, E_i)|^2 \leq \frac{1}{2}$}

\noindent In this case, $ \rho \geq \frac{1}{2} $. Combine inequality (\ref{eq:system2-final}) with (\ref{eq:fi}) and (\ref{eq:fihi^2}), we obtain:

	\begin{align}\label{eq:case1main}
		F^{ij} P_{ij} &\geq -\frac{16F}{\kappa_n} - C \left(1 + \beta + \alpha \right) 
		+ \frac{F h_{11}}{\log h_{11}} + \frac{C \alpha}{\kappa_n}.
	\end{align}
	
	Choosing $ \alpha $ sufficiently large, we obtain:
	\begin{align}
		F^{ij} P_{ij} &\geq - C \left(1 + \beta + \alpha \right) + \frac{F h_{11}}{\log h_{11}}.
	\end{align}
	
This provides the estimate for this case. Hence, it remains to discuss $\sum_{i=1}^n |(X, E_i)|^2 > \frac{1}{2}$.

~\\
\noindent\textbf{\textit{Case 2: }$\sum_{i=1}^n |(X, E_i)|^2 > \frac{1}{2}$ and suppose $ \sqrt{\log \kappa_1} \geq C \kappa_{n-1}$}

\noindent
    By lemma \ref{lm:lambda_nbound}, we have:
	\begin{align*}
		\left( \frac{1}{\kappa_{n-1}} + \dots + \frac{1}{\kappa_1} \right) \leq \dfrac{\kappa_n}{F} \leq C(n) \left( \frac{1}{\kappa_{n-1}} + \dots + \frac{1}{\kappa_1} \right)
	\end{align*}
    Using the condition of case 2, we obtain:
	\begin{align}
		\dfrac{\kappa_n}{F} \geq \frac{1}{\kappa_{n-1}} + \dots + \frac{1}{\kappa_1} \geq \dfrac{1}{\kappa_{n-1}} \geq \frac{C}{\sqrt{\log h_{11}}}, \label{eq:case3.1solu}
	\end{align}
     Substituting relation \eqref{eq:fi}, (\ref{eq:relaofrhoandlambda_1}) and (\ref{eq:case3.1solu}) into (\ref{eq:system2-final}), we derive:

\begin{align*}
	F^{ij} P_{ij}
&\geq -\dfrac{4 \sum_{i} F^{ii}}{\rho} - C(\frac{1}{\rho} + \beta + \alpha ) +\frac{F h_{11}}{\log h_{11}} +\alpha \frac{F^{ii} h_{ii}^2}{(\nu, E_{n+1})^2} \notag \\
&\geq -\dfrac{8 F}{\kappa_n\rho} - C(\frac{1}{\rho} + \beta + \alpha ) +\frac{F h_{11}}{\log h_{11}} +\alpha \frac{F^{ii} h_{ii}^2}{(\nu, E_{n+1})^2} \\
&\geq -C\log {\kappa_1} - C(\sqrt{\log {\kappa_1}} + \beta + \alpha ) +\frac{F h_{11}}{\log h_{11}} +\alpha \frac{F^{ii} h_{ii}^2}{(\nu, E_{n+1})^2}
\end{align*}

For sufficiently large $ \kappa_1 $, choosing $ \frac{1}{2} F \kappa_1 > C \log^2 \kappa_1 $, gives the estimate for case 2.

~\\
We know from inequality (\ref{eq:moequvalent-b}) that it is impossible that all $|b_i|$ are very small. In case 3, we will divide it to two sub-cases by whether there is a large enough $|b_j|$ for $\{j=1,2,\cdots,n-1\}$ or not.

\noindent \textit{\textbf{Case 3.1}}: \textbf{Suppose $ \sqrt{\log \kappa_1} < C \kappa_{n-1} $, and $\sum_{i=1}^n |(X, E_i)|^2 > \frac{1}{2}$. Let $d>0$ be a small constant to be determined later and only $|b_n|>d$ among all $|b_i|$.}

Recall that $\textbf{b}=T\textbf{a}$, the matrix $T$ is symmetric, positive-definite and invertible. Hence we have $\textbf{a}=T^{-1}\textbf{b},$ and by the Sherman--Morrison formula,
\[
T^{-1}
=
\frac{1}{(\nu,E_{n+1})}
\left(I-\textbf{c}\textbf{c}^t\right).
\]

Hence we can calculate
\begin{align}
|a_n|
&=
\left|
\frac{1}{(\nu,E_{n+1})}
\left[
b_n
-
(e_n,E_{n+1})
\sum_{j=1}^{n}
b_j(e_j,E_{n+1})
\right]
\right|
\notag\\
&=
\left|
\frac{1-(e_n,E_{n+1})^2}{(\nu,E_{n+1})}b_n
-
\frac{(e_n,E_{n+1})}{(\nu,E_{n+1})}
\sum_{j=1}^{n-1}
b_j(e_j,E_{n+1})
\right|
\notag\\
&=
\left|
\frac{(\nu,E_{n+1})^2+\displaystyle\sum_{j=1}^{n-1}(e_j,E_{n+1})^2}
{(\nu,E_{n+1})}
\,b_n
-
\frac{(e_n,E_{n+1})}{(\nu,E_{n+1})}
\sum_{j=1}^{n-1}
b_j(e_j,E_{n+1})
\right|
\notag\\
&\ge
(\nu,E_{n+1})|b_n|
-
\left|
\frac{(e_n,E_{n+1})}{(\nu,E_{n+1})}
\sum_{j=1}^{n-1}
b_j(e_j,E_{n+1})
\right|.
\notag
\end{align}

We now choose $d>0$ sufficiently small so that
\begin{align}
\frac{1}{2}(\nu,E_{n+1})|b_n|
&\ge
\left|
\frac{(e_n,E_{n+1})}{(\nu,E_{n+1})}
\sum_{j=1}^{n-1}
b_j(e_j,E_{n+1})
\right|.
\notag
\end{align}

A sufficient condition is
\begin{align}\label{eq:rangeofd}
d
\le
\frac{(\nu,E_{n+1})^2}
{\sqrt{2C_7(n-1)\bigl[(\nu,E_{n+1})^4+4\bigr]}}.
\end{align}

Consequently, if \eqref {eq:rangeofd} is satisfied,
\begin{align}\label{eq:anlower}
|a_n|
\ge
\frac{1}{2}(\nu,E_{n+1})|b_n|.
\end{align}

Combining \eqref{eq:anlower} with \eqref{eq:lambda_n-1,nbound+} and \eqref{eq:fi_ele}, we estimate a positive term in \eqref{eq:system2-final} as follows
\begin{align}\label{eq:case3.1main}
\hat{c}(n)\log h_{11}\cdot
\frac{4\sum_{i=1}^{n}F^{ii}a_i^2}{\rho^2}
&\ge
\hat{c}(n)\log h_{11}\cdot
\frac{4F^{nn}a_n^2}{\rho^2}
\notag\\
&\ge
\hat{c}(n)\log h_{11}\cdot
\frac{CFd^2}{\kappa_n\rho^2}.
\end{align}

Substituting \eqref{eq:case3.1main} into \eqref{eq:system2-final} and using $\rho\le\frac12$ together with \eqref{eq:fi}, we obtain

\begin{align}
0
\ge&
-\frac{4\sum_iF^{ii}}{\rho}
-C\left(\frac{1}{\rho}+\beta+\alpha\right)
+\hat{c}(n)\log h_{11}\cdot
\frac{4F^{ii}a_i^2}{\rho^2}
\notag\\
&
+\frac{\hat{c}(n)\log h_{11}}{2}\cdot
F^{ii}h_{ii}^2
\left(
\beta\frac{b_i}{(\nu,E_{n+1})^2}
+
2\alpha\frac{(e_i,E_{n+1})}{(\nu,E_{n+1})^3}
\right)^2
\notag\\
&
+\frac{Fh_{11}}{\log h_{11}}
+\alpha\frac{F^{ii}h_{ii}^2}{(\nu,E_{n+1})^2}
\notag\\
\ge&
-\frac{8F}{\kappa_n\rho}
-C\left(\frac{1}{\rho}+\beta+\alpha\right)
+\hat{c}(n)\log h_{11}\cdot
\frac{CFd^2}{\kappa_n\rho^2}
\notag\\
&
+\frac{Fh_{11}}{\log h_{11}}
+\alpha\frac{F^{ii}h_{ii}^2}{(\nu,E_{n+1})^2}.
\notag
\end{align}

Since $h_{11}$ is sufficiently large and $\rho\le\frac12$, we obtain
\begin{align}
0
\ge&
-C(1+\beta+\alpha)
+\frac{Fh_{11}}{\log h_{11}}
+\alpha\frac{F^{ii}h_{ii}^2}{(\nu,E_{n+1})^2}.
\notag
\end{align}

This completes the proof of Case~3.1.

~\\
\noindent \textit{\textbf{Case 3.2}}: \textbf{Suppose $ \sqrt{\log \kappa_1} < C \kappa_{n-1} $, and $\sum_{i=1}^n |(X, E_i)|^2 > \frac{1}{2}$. Assume there exists a $j\in \{1,2,\cdots,n-1\}$ such that $|b_j|>d$, where the small constant $d$ is the same as that in case 3.1.}

From inequality \eqref{eq:rangeofd}, there exists such a constant $d>0$, independent of $\alpha$ and $\beta$.

Choose $\alpha$ sufficiently large so that all the previous conditions are satisfied. Then choose $\beta$ so that
\begin{align}\label{eq:betalower}
\beta d
\ge
\frac{4\alpha}{(\nu,E_{n+1})}.
\end{align}

Together with inequality \eqref{eq:fi_ele}, we estimate another positive term in \eqref{eq:system2-final} as follows, the $j$ here is the fixed index mentioned in the assumption of case~3.2:
\begin{align}\label{eq:case3.2main}
&\frac{\hat{c}(n)\log h_{11}}{2}\cdot
\sum_{i=1}^{n}F^{ii}h_{ii}^{2}
\left(
\beta\frac{b_i}{(\nu,E_{n+1})^{2}}
+
2\alpha\frac{(e_i,E_{n+1})}{(\nu,E_{n+1})^{3}}
\right)^{2}
\notag\\
\ge&
\frac{\hat{c}(n)\log h_{11}}{2}\cdot
F^{jj}h_{jj}^{2}
\left(
\beta\frac{b_j}{(\nu,E_{n+1})^{2}}
+
2\alpha\frac{(e_j,E_{n+1})}{(\nu,E_{n+1})^{3}}
\right)^{2}
\notag\\
\ge&
\frac{\hat{c}(n)\log h_{11}}{2}
\cdot
\frac{F^{2}h_{jj}^{2}}
{\kappa_n\kappa_j^{2}}
\left(
\frac{\beta d}
{2(\nu,E_{n+1})^{2}}
\right)^{2}
\notag\\
\ge&
\frac{\hat{c}(n)\beta^{2}d^{2}F^{2}}
{8(\nu,E_{n+1})^{4}}
\cdot
\frac{\log h_{11}}
{\kappa_n}.
\end{align}

Substituting \eqref{eq:relaofrhoandlambda_1} and \eqref{eq:case3.2main} into \eqref{eq:system2-final}, together with \eqref{eq:fi} again, we obtain
\begin{align}
0
\ge&
-\frac{4\sum_iF^{ii}}{\rho}
-C\left(\frac{1}{\rho}+\beta+\alpha\right)
+\hat{c}(n)\log h_{11}\cdot
\frac{4F^{ii}a_i^2}{\rho^2}
\notag\\
&
+\frac{\hat{c}(n)\log h_{11}}{2}\cdot
F^{ii}h_{ii}^{2}
\left(
\beta\frac{b_i}{(\nu,E_{n+1})^{2}}
+
2\alpha\frac{(e_i,E_{n+1})}
{(\nu,E_{n+1})^{3}}
\right)^{2}
\notag\\
&
+\frac{Fh_{11}}{\log h_{11}}
+\alpha\frac{F^{ii}h_{ii}^{2}}
{(\nu,E_{n+1})^{2}}
\notag\\
\ge&
-\frac{8CF\sqrt{\log h_{11}}}
{\kappa_n}
-C\left(\sqrt{\log h_{11}}+\beta+\alpha\right)
+\frac{\hat{c}(n)\beta^{2}d^{2}F^{2}}
{8(\nu,E_{n+1})^{2}}
\cdot
\frac{\log h_{11}}
{\kappa_n}
\notag\\
&
+\frac{Fh_{11}}{\log h_{11}}
+\alpha\frac{F^{ii}h_{ii}^{2}}
{(\nu,E_{n+1})^{2}}.
\notag
\end{align}

Since $h_{11}$ is sufficiently large and both $\beta$ and $d$ are positive constants, we derive
\begin{align}
0
\ge&
-C(1+\beta+\alpha)
+\frac{Fh_{11}}
{\log h_{11}}
+\alpha\frac{F^{ii}h_{ii}^{2}}
{(\nu,E_{n+1})^{2}}.
\notag
\end{align}

This completes the proof for case~3.2.

~\\
We finish the proof by summarizing the order of parameter selection. We first choose $d>0$ sufficiently small to satisfy \eqref{eq:rangeofd}. Then, we choose $\alpha$ sufficiently large in \eqref{eq:system2-3} and \eqref{eq:case1main} to satisfy the following steps, and finally choose $\beta$ sufficiently large to satisfy \eqref{eq:betalower}. With these choices, all the parameters in the test function $P(X)$ are fixed uniformly, and the proof of Theorem~\ref{thm:1.1} is complete.

\end{proof}

\section*{Acknowledgements}
This work was carried out while the author was a Master's student at McMaster University, and constitutes the author's Master's thesis. The current version also corrects certain errors present in the original thesis manuscript. The author would like to thank Professor Siyuan Lu for helpful discussions and for bringing many related works to our attention. He is also deeply grateful to Bin Wang for numerous important discussions and for his careful suggestions regarding the correctness of this work. In particular, the author sincerely thanks Professor Guohuan Qiu for kindly reviewing the manuscript with great care.

\bibliography{main}
\bibliographystyle{amsplain}
\end{document}